\documentclass[12pt]{article}
\usepackage{amsmath,amssymb,amsthm}
\usepackage{bbm}
\usepackage{url}
\usepackage[arrow,matrix,curve]{xy}
\usepackage{tikz}
\usepackage{tikz-cd}
\usepackage{hyperref}

\newtheorem*{klem}{Key-Lemma}
\newtheorem{prop}{Proposition}
\newtheorem{cor}{Corollary}
\newtheorem*{thm*}{Theorem}
\newtheorem*{prop*}{Proposition}
\newtheorem*{cor*}{Corollary}
\theoremstyle{definition}
\theoremstyle{definition}
\newtheorem{rmk}{Remark}
\theoremstyle{definition}

\newcommand{\ZZ}{\mathbb Z}

\newcommand{\RR}{\mathbb R}
\newcommand{\QQ}{\mathbb Q}
\newcommand{\CC}{\mathbb C}

\def\diag{\mathop{\mathrm{diag}}\nolimits}

\renewcommand{\Re}{{\rm Re}}
\renewcommand{\Im}{{\rm Im}}

\begin{document}

\author{Hanno von Bodecker\footnote{Fakult{\"a}t f{\"u}r Mathematik, Universit{\"a}t Bielefeld, Germany}  \and  Sebastian Thyssen\footnote{Fakult\"at f\"ur Mathematik, Ruhr-Universit\"at Bochum, Germany} }

\title{On $p$--local Topological Automorphic Forms for $U(1,1;\mathbb{Z}[i])$ }

\date{}

\maketitle

\begin{abstract}
We present a new flavor of  TAF-type (co)homology theories, which are $p$--local of height two and based on the isometry group of the odd unimodular hermitian lattice of signature $(1,1)$ over the Gaussian integers. Using a suitable family of hyperelliptic curves, we explicitly construct a genus to automorphic forms, prove an integrality statement and verify Landweber's criterion.
\end{abstract}

\section{Introduction}

In \cite{Landweber:1995sw} Landweber, Ravenel and Stong used bordism theory to produce complex oriented cohomology theories arising from formal groups of elliptic curves. These so called \textit{elliptic} cohomology theories were a first outcome of the analysis of the connection between stable homotopy theory and algebraic geometry, in particular the study of abelian varieties, via the theory of $1$-dimensional formal groups. 

The natural way to associate formal groups to abelian schemes is by completing the abelian schemes against the zero section, the formal group of an $n$-dimensional abelian scheme naturally being $n$-dimensional. Elliptic curves thus provide an obvious starting point in the investigation of the above connection.

The study of elliptic cohomology culminated in the construction of the spectrum $Tmf$ of \textit{topological modular forms} in its various flavors, which, roughly speaking, realizes the canonical $1$-dimensional formal group on the (compactified) moduli stack of elliptic curves. Though not elliptic, the theory is commonly understood as `universal elliptic' theory.

\bigskip

The investigation of the connection between stable homotopy theory and abelian varieties was continued  in \cite{Behrens:2010aa} by Behrens and Lawson via the theory of \textit{topological automorphic forms}. This theory is a generalization of the theory of topological modular forms in the sense of classifying higher dimensional abelian schemes via unitary Shimura stacks. 

The obvious problem of how to associate a $1$-dimensional formal group to such Shimura stack is being solved by equipping the abelian schemes in  question with enough extra structure to split the associated formal groups and ensure that one summand is of dimension one. 
The construction of $TAF$-spectra is similar to the construction of $TMF$, with Behrens and Lawson using a realization result of Lurie's to define a sheaf of $E_\infty$-ring spectra on $p$-completions of the particular $U(n-1,1)$-Shimura varieties, the spectrum of topological automorphic forms being the global sections object. 

\bigskip

The purpose of the present report is to give an {explicit} geometric construction of the $p$-local theory of topological automorphic forms for $U(1,1;\mathbb{Z}[i])$, the isometry group of the odd unimodular lattice of signature $(1,1)$ over the Gaussian integers. 
Following \cite{Landweber:1995sw} we present a suitable genus in terms of Legendre polynomials and use bordism theory to produce $TAF$-type cohomology theories.

\bigskip

We start by collecting relevant well-known material about symplectic and unitary groups and abelian varieties parametrized by quotients of such. 

Turning to the Shimura variety, we continue with a careful investigation of the geometry of the complex points. We determine a fundamental domain for the unitary action on the universal cover, and then present a twisted version of the standard embedding of unitary into symplectic groups, which enables us to pin down the splitting of the Lie algebra of the parametrized abelian varieties. The explicit calculation of the ring of automorphic forms puts us in the position to completely describe the complex geometry.

Our strategy is to access the (integral model of the) Shimura variety via a certain family of hyperelliptic curves of genus two with an action of the fourth roots of unity. The action distinguishes a differential on the curve, thus splits the cotangent bundle, and extends to an action of the Gaussian integers on the Jacobian that is subject to the required $(1,1)$-signature condition. 

Over the complex numbers our family of curves is the universal family of such curves. Having determined the integral subring of the ring of automorphic forms, we can describe the family in terms of automorphic forms.

\bigskip

We address the integrality of the involved formal groups. Following Behrens and Lawson, we argue that, for split prime $p$ in the Gaussian integers, the formal group associated to the abelian two-folds splits over the $p$-completion of the Shimura variety into two $1$-dimensional formal summands, each defined over the $p$-adic integers.  

Furthermore, as we deal with two-folds coming from curves, we obtain a rational formal group, with logarithm given by integrating the distinguished differential. 

Choosing a suitable local parameter at the point at infinity on our (family of) curve(s) allows us to pull back the $1$-dimensional formal summand of the formal group of the Jacobian and thus describe it via the curve differential, which in turn, using the explicit presentation of the curve, can be described in terms of (Legendre polynomials $P_k$ in) automorphic forms.

\begin{thm*} 
Let $M^{R}_*(G)  $ denote the ring of automorphic forms for $U(1,1;\mathbb{Z}[i])$ with expansions in $R[\![q^{1/2}]\!]$ and let us define a rational genus
\[
\varphi^{L}\colon MU_{*}^{\QQ}\to M_{*}^{\QQ}(G)
\]
by means of 
\begin{align*} \log_{\varphi^{L}}'(x)
= \sum_{k=0}^\infty \varphi^{L}[\CC P^{k}]\, x^{k} 
= \sum_{k=0}^\infty P_{k}(\alpha,\beta)\, x^{4k} .
\end{align*}
Then for each prime $p\equiv1\!\mod4$,  the image of $BP_*$ under $\varphi^L$ is contained in the subring of $p$-local automorphic forms:
\[
\varphi^L(BP_*) \subset M^{\mathbb{Z}_{(p)}}_*(G)  \cong\mathbb{Z}_{(p)}[\alpha, \beta].
\]
\end{thm*}

Note that for $n=2$ the Torelli map, which on points embeds the curve into its Jacobian, 
is well behaved, i.e.\ an immersion, on the level of stacks. On coarse moduli schemes the image of our family even provides a dense sublocus in the Shimura variety. More precisely, for the odd lattice there is only one point where the abelian two-fold is not the Jacobian of a hyperelliptic curve. 

It is a consequence of our investigation of complex points that in this point (corresponding to $\beta=0$) the two-fold is the product of two elliptic curves (each) with Gaussian action. Further, in this point the group law associated to $\varphi^{L}$ reduces to (a variant of) Euler's formal group law $F_{E}$,
\[
F_{E}(x,y)=\frac{x\sqrt{1-2\alpha y^{4}}+y\sqrt{1-2\alpha x^{4}}}{1+2\alpha x^{2}y^{2}}.
\]

\bigskip

Finally, we show that our genera take values in Landweber exact $BP_*$-algebras and that we can thus produce $TAF$-type (co)homology theories.

\begin{cor*}
Let $p=5$. Then $\varphi^{L}$ gives $M_{*}^{\ZZ_{(5)}}(G)[\Delta_{G}^{-1}]$ the structure of a Landweber exact $BP_*$-algebra of height two. In particular, the functors
\[
TAF^{U(1,1;\ZZ[i])}_{(5),*}(\cdot)=BP_{*}(\cdot)\otimes_{\varphi^{L}}M^{\ZZ_{(5)}}_{*}(G)[\Delta_{G}^{-1}]
\]
define a homology theory.
\end{cor*}

This is also true for $p=13$ and conjecturally for all $p\equiv 1\!\mod 4$. Anyhow, we content ourselves with a slightly weaker statement over a marginally larger ring:
\begin{cor*}
Let $p\equiv5\mod8$, then $v_{2}$ is non-zero mod $(p,v_{1})$. In particular, $M_{*}^{\ZZ_{(p)}}(G)[v_{2}^{-1}]$ admits the structure of a Landweber exact algebra of height two.
\end{cor*}

\subsubsection*{Remarks}

There has been an analysis of the \textit{even} unimodular lattice of signature $(1,1)$ over the Gaussian integers in \cite{Behrens:2011aa}. 

Note further, that our methods work analogously for the unique unimodular lattice over the Eisenstein integers. For a more detailed account we refer the reader to the forthcoming paper \cite{TAFviaCurves}.

\subsubsection*{Acknowledgements}

The second author thanks the DFG SPP 1786 for financial support.


\section{Preliminaries}	
\label{SectionPrelim}
In this section we collect the relevant material for the results we produce in the next section. We start giving a brief recollection of the basic unitary and symplectic groups involved, and explain how the former can be embedded into the latter in a way compatible with the abelian varieties parametrized by quotients of these groups.
\subsection{Symplectic groups and abelian varieties} 
The matrix $\left( \begin{smallmatrix} &  \mathbbm{1}\\ -\mathbbm{1} & \end{smallmatrix} \right)$, where $\mathbbm{1}$ denotes the unit matrix of rank $n$,
defines a skew bilinear form on $\RR^{2n}$. The symplectic group is defined as
\[
Sp(n)=\{g\in GL(2n;\RR):g^{tr}\left( \begin{smallmatrix} &  \mathbbm{1}\\ -\mathbbm{1} & \end{smallmatrix} \right) g=\left( \begin{smallmatrix} &  \mathbbm{1}\\ -\mathbbm{1} & \end{smallmatrix} \right)\};
\]
it is a non-compact real Lie group of dimension $2n^{2}+n$ and real rank $n$. Moreover, we have the inversion formula
\[
\left(\begin{smallmatrix}A&B\\C&D\end{smallmatrix}\right)^{-1}=\left(\begin{smallmatrix}D^{tr}&-B^{tr}\\-C^{tr}&A^{tr}\end{smallmatrix}\right).
\]
Since $\left( \begin{smallmatrix} &  \mathbbm{1}\\ -\mathbbm{1} & \end{smallmatrix} \right)^{-1}=-\left( \begin{smallmatrix} &  \mathbbm{1}\\ -\mathbbm{1} & \end{smallmatrix} \right)$, the group $Sp(n)$ is closed under transposition. Thus, $g\mapsto (g^{tr})^{-1}$ is an involutive automorphism. Let $$Sp(n)\cap O(2n)\cong U(n)$$ be the fixed point set; it is a maximal compact subgroup, and there is a well-known identification between the associated symmetric space and the Siegel space $\mathbb{S}_{n}$ consisting of symmetric complex matrices with positive-definite imaginary part,
\[
Sp(n)/U(n)\cong \mathbb{S}_{n}=\{\Omega=\Re(\Omega)+i\Im(\Omega):\Omega^{tr}=\Omega,\ \Im(\Omega)>0\}.
\]
The Siegel upper half space $\mathbb{S}_n$ is a hermitian symmetric domain well known to parametrize (complex) abelian varieties equipped with a principal polarization: For $\Omega\in\mathbb{S}_{n}$, the columns of the matrix $( \mathbbm{1},\Omega)$ constitute a $\ZZ$--module basis for a symplectic lattice 
\[
\Lambda = \mathbb{Z}^n \oplus \Omega\mathbb{Z}^n
\]
in $\mathbb{C}^n$, and the quotient $\CC^{n}/\Lambda$ becomes a principally polarized abelian variety.  Letting 
\[
\sigma\colon Sp(n)\to Sp(n), \quad \sigma(\begin{smallmatrix}A&B\\C&D\end{smallmatrix})=(\begin{smallmatrix}A&-B\\-C&D\end{smallmatrix})
\]
denote the automorphism introducing a sign in the off-diagonal blocks, we have the left action
\[
( \mathbbm{1},\Omega)\sigma(\begin{smallmatrix}A&B\\C&D\end{smallmatrix})^{-1}=(D^{tr}+\Omega C^{tr}, B^{tr}+\Omega A^{tr}),
\]
allowing us to recover the usual fractional linear transformation on $\mathbb{S}_{n}$,
\[
\left(\begin{smallmatrix}A&B\\ C&D\end{smallmatrix}\right)\cdot\Omega=(A\Omega+B)(C\Omega+D)^{-1}.
\]
Clearly, $Sp(n;\ZZ)$--equivalent points in Siegel space determine isomorphic abelian varieties.
Finally, we remark that there is a left action of $Sp(n;\ZZ)$ on $\mathbb{S}_{n}\times\CC^{n}$ given by
\begin{equation}\label{almost a bundle of Lie algebras}
\left(\begin{smallmatrix}A&B\\ C&D\end{smallmatrix}\right)\cdot(\Omega,\vec{w})=\big((A\Omega+B)(C\Omega+D)^{-1},((C\Omega+D)^{tr})^{-1}\vec{w} \big);
\end{equation}
the quotient tries to be the bundle of Lie algebras of the respective abelian varieties.
\subsection{The unitary group $U(1,1;\ZZ[i])$}
There are two isometry classes of unimodular lattices of signature $(1,1)$ over the Gaussian integers $\ZZ[i]$. Here, we focus on the odd lattice, represented by the diagonal Gram matrix $H=\left(\begin{smallmatrix}1&\\&-1\end{smallmatrix}\right)$. Its isometry group is 
\[
U(1,1;\ZZ[i])\subset U(1,1)= \{g\in GL_2(\mathbb{C}) : \bar{g}^{tr}Hg = H \}
\]
and forms a cofinite-volume (but not cocompact) lattice in the reductive Lie group $U(1,1)$. A maximal compact subgroup of this Lie group is given by
\[
U(1,1)\cap U(2)\cong U(1)\times U(1).
\]
Clearly, this group constitutes the fixed point set of the involutive automorphism $g\mapsto (\bar g^{tr})^{-1}=HgH$. The associated symmetric space admits a bounded realization as the open ball $\mathbb{B}_1$ in $\CC$,
\[
U(1,1)/(U(1,1)\cap U(2))\cong\mathbb{B}_{1}=\{z\in\CC:\bar{z}^{tr}z<1\}.
\]

More conceptually, $\mathbb{B}_{1}$ parametrizes the lines in $\CC^{2}$ which are negative-definite w.r.t.\ $H$ (i.e.\ an open subset in $\CC P^{1}$, $z\mapsto\left[\begin{smallmatrix}z\\1\end{smallmatrix}\right])$.
\subsection{Standard embedding}
Recall that $U(1,1;\ZZ[i])$ is the isometry group of a unimodular hermitian Gaussian lattice. The underlying $\ZZ$--lattice carries a unimodular skew symmetric $\ZZ$--bilinear form, $(x,y)\mapsto\Im(\bar{x}^{tr}Hy)$. Thus, $U(1,1;\ZZ[i])$ can be embedded into $Sp(2;\ZZ)$. Specifically, we choose the embedding $$\rho\colon U(1,1;\ZZ[i])\to Sp(2;\ZZ),$$
given by
\[
g=\Re(g)+i\Im(g)\mapsto\begin{pmatrix}\Re(g)&\Im(g)H\\-H\Im(g)&H\Re(g)H\end{pmatrix},
\]
which induces the (holomorphic) map $\rho_{*}\colon\mathbb{B}_{1}\to\mathbb{S}_{2}$ given by
\[
z\mapsto\Omega_{z}=\frac{i}{1-z^{2}}\begin{pmatrix}1+z^{2}&2z\\2z&1+z^{2}\end{pmatrix}.
\]
To see this, write $z=e^{i\phi}\tanh\psi$ and consider the image of
\[
g=\begin{pmatrix}\cosh\psi&e^{i\phi}\sinh\psi\\e^{-i\phi}\sinh\psi&\cosh\psi\end{pmatrix}\in SU(1,1),
\]
letting it act on $\diag(i,i)\in\mathbb{S}_{2}$.

To characterize the image, let $J=(\begin{smallmatrix}&-H\\H&\end{smallmatrix})$; then it is straightforward to check that the homomorphism $\rho$ induces an identification
\begin{equation}\label{U in Sp}
U(1,1;\ZZ[i])\cong\{g\in Sp(2;\ZZ):\sigma(g)J=J\sigma(g)\}.
\end{equation}
Clearly, $J$ yields an automorphism of the (symplectic) lattice spanned by the columns of $( \mathbbm{1},i \mathbbm{1})$, thereby providing a module structure over the Gaussian integers $\ZZ[i]$. Extending \eqref{U in Sp} to a homomorphism of real Lie groups, it follows that the same holds true for all the lattices in its $\rho(U(1,1))$--orbit.

Furthermore, the action of $J$ can also be described using left multiplication by  a complex two-by-two matrix:
\begin{equation}\label{lifted automorphism}
( \mathbbm{1},\Omega_{z})\cdot J=(\Omega_{z}H,-H)=iR_{z}\cdot( \mathbbm{1},\Omega_{z}),
\end{equation}
where the matrix
\[
R_{z}=-i\Omega_{z}H=\frac{1}{1-z^{2}}\begin{pmatrix}1+z^{2}&-2z\\2z&-1-z^{2}\end{pmatrix}
\]
describes an involution of determinant $-1$; we remark that
\[
R_{z}\cdot\begin{pmatrix}z\\1\end{pmatrix}=-\begin{pmatrix}z\\1\end{pmatrix}.
\]
From \eqref{lifted automorphism} it follows that the action of $iR_z$ on $\CC^{2}$ descends to an automorphism of order four on the abelian variety defined by $\Omega_z$. This makes clear that the abelian varieties parametrized by the sublocus $\rho_*(\mathbb{B}_1)\subset \mathbb{S}_2$ carry the additional structure of a $\mathbb{Z}[i]$-action.

\begin{rmk}
In fact, the group of automorphisms of these abelian varieties is larger: Let $H'=\left(\begin{smallmatrix}&1\\1&\end{smallmatrix}\right)$. Then the involution $S'=\left(\begin{smallmatrix}H'&\\&H'\end{smallmatrix}\right)$ lies in $Sp(2;\ZZ)$ (although it does {\em not} lie in the image of $\rho$), stabilizes each element of $\rho_{*}(\mathbb{B}_1)$ and induces an automorphism on the corresponding varieties. 

Due to $S'JS'=-J=J^{-1}$, the automorphism group therefore contains a dihedral group with eight elements. We remark that since $S'$ is not in the image of $\rho$, it would be better to work with a group larger than $U(1,1)$, namely the group $GU(1,1)$ preserving the hermitian form $H$ only up to scale, cf. \cite{Shimura:1963aa}.
\end{rmk}
\section{The geometry of the complex Shimura variety} \label{SectionGeomShimVar}
We analyze the complex Shimura variety in more detail:
Applying a Cayley transform to the unit ball to work with the more familiar upper half plane, we determine a fundamental domain. Giving a twisted version of the standard embedding, we make manifest the splitting of the Lie algebra of the associated abelian varieties. We determine the ring of automorphic forms and use this to draw some conclusions about the geometry of the quotient space.

\subsection{Cayley transform}

The common analytic picture interprets modular forms as functions on the upper half space. In order to import this familiar point of view into our theory, we apply the Cayley transform to the unit ball. 
\bigskip

The (left) action of $g_{0}=\frac{1}{\sqrt2}\left(\begin{smallmatrix}1&i\\i&1\end{smallmatrix}\right)\in SU(2)$ on $\CC P^{1}$ induces the Cayley transform $\mathbb{B}_{1}\to\mathbb{H}_1=\{\tau\in\CC:\Im(\tau)>0\}$,
\[
z\mapsto \tau=\frac{z+i}{1+iz}=i\frac{1-iz}{1+iz}.
\]
On the level groups, the transform takes the form $U(1,1)\to g_{0}U(1,1)g_{0}^{-1}$,
\[
\begin{pmatrix}a&b\\c&d\end{pmatrix}\mapsto\frac{1}{2}\begin{pmatrix}a+ic-ib+d&-ia+c+b+id\\ia+c+b-id&a-ic+ib+d\end{pmatrix}.
\]

Clearly, the surjective homomorphism $\det\colon U(1,1;\ZZ[i])\to(\ZZ[i])^{\times}$ admits a section; thus, we have a decomposition
\[
U(1,1;\ZZ[i])\cong SU(1,1;\ZZ[i])\rtimes C_{4},
\]
where $C_{4}$ denotes a cyclic group of order four. Furthermore, we have:
\begin{prop} The Cayley transform induces an isomorphism
\[
SU(1,1;\ZZ[i])\cong\Gamma_{\theta}=\{\left(\begin{smallmatrix}a&b\\c&d\end{smallmatrix}\right)\in SL(2;\ZZ):ab\equiv cd\equiv0(2)\}
\]
\end{prop}
\begin{proof} 
In $U(1,1)$ we have
\[
\begin{pmatrix}a&b\\c&d\end{pmatrix}^{-1}=\begin{pmatrix}\bar a&-\bar c\\-\bar b&\bar d\end{pmatrix}
\]
while in $SL(2;\CC)$, we have
\[
\begin{pmatrix}a&b\\c&d\end{pmatrix}^{-1}=\begin{pmatrix}d&-b\\-c&a\end{pmatrix}
\]
Thus, working in $SU(1,1)$, we may put $d=\bar a$, $b=\bar c$, and therefore
\[
g_{0}\begin{pmatrix}a&\bar c\\c&\bar a\end{pmatrix}g_{0}^{-1}=\frac{1}{2}\begin{pmatrix}a+ic-i\bar c+\bar a&-ia+c+\bar c+i\bar a\\ia+c+\bar c-i\bar a&a-ic+i\bar c+\bar a\end{pmatrix}\in SL(2;\RR).
\]
Furthermore, we know that $|a|^{2}-|c|^{2}=1$; thus, if $a,c\in\ZZ[i]$, then the rational integers $\Re(c)+\Im(a)$ and $\Re(a)+\Im(c)$ have different parity, establishing injectivity. Surjectivity follows from the fact that
\begin{equation*}
g_{0}\begin{pmatrix}1+i&1\\1&1-i\end{pmatrix}g_{0}^{-1}=\begin{pmatrix}1&2\\0&1\end{pmatrix},\ g_{0}\begin{pmatrix}i&\\&-i\end{pmatrix}g_{0}^{-1}=\begin{pmatrix}&1\\-1&\end{pmatrix},
\end{equation*}
are well-known  generators for the  group $\Gamma_{\theta}\subset SL(2;\ZZ)$ (cf., e.g., \cite{Freitag:2009aa}).
\end{proof}

Recall that $\Gamma_{\theta}$ is a subgroup of index three in $SL(2;\ZZ)$ and that a fundamental domain for the action of $\Gamma_{\theta}$ on the upper half-plane $\mathbb{H}_1$ can be described as the vertical strip between $-1$ and $1$ of elements $\tau\in\mathbb{H}_1$ of (euclidean) norm $|\tau|\geq1$. Choosing $$g_0 \left( \begin{smallmatrix} i & \\ & 1  \end{smallmatrix}\right)g_0^{-1} = \tfrac{1+i}{2}\left(\begin{smallmatrix}    1 &1 \\ -1 & 1\end{smallmatrix}\right)$$ as generator of the (transformed) $C_4$, which induces the map
\begin{equation}\label{additional transformation}
\tau\mapsto\frac{1+\tau}{1-\tau}=-\frac{2}{\tau-1}-1,
\end{equation}
it is easy to see that a fundamental domain for the Cayley transform of $U(1,1;\mathbb{Z}[i])$ is given by the part of a vertical strip $-1\leq x\leq 1$ lying above the two circles of radius $\sqrt{2}$ around $-1$ and $1$.
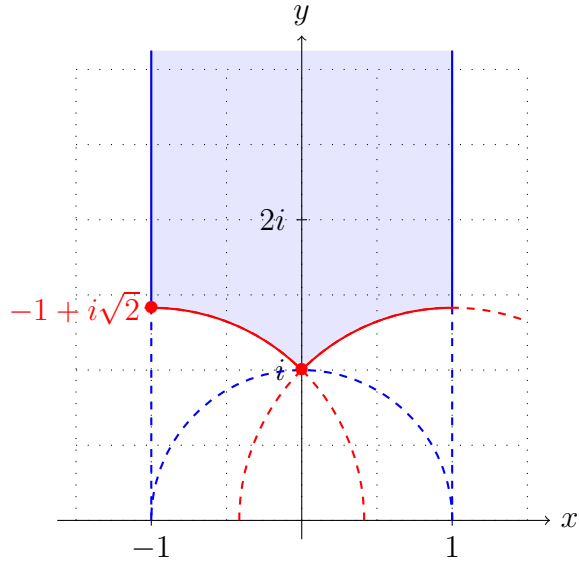
\begin{center}
\begin{tikzpicture}[scale=1.0]
\filldraw[fill=blue!10!white, fill opacity=20] 
 (-2,6.25) -- (-2,2.83) -- (-2,2.83) arc(90:45:2.83) -- (0,2) arc (135:90:2.83) -- (2,2.83) -- (2,6.25) ;
\draw[loosely dotted] (-3,0) grid (3,6);
\draw[->] (-3.25,0) -- (3.3,0) node[right] {$x$};
\foreach \x/\xtext in {-2/-1,  2/1}
\draw[shift={(\x,0)}] (0pt,2pt) -- (0pt,-2pt) node[below] {$\xtext$};
\draw[->] (0,-0.25) -- (0,6.45) node[above] {$y$};
\foreach \y/\ytext in { 2/i, 4/2i}
    \draw[shift={(0,\y)}] (2pt,0pt) -- (-2pt,0pt) node[left] {$\ytext$};
\draw[thick, blue] (-2,2.83) -- (-2,6.25);
\draw[thick, blue, dashed] (-2,0) -- (-2,2.83);
\draw[thick, blue] (2,2.83) -- (2,6.25);
\draw[thick, blue, dashed] (2,0) -- (2,2.83);
\draw[blue,thick, dashed]      (2,0) arc (0:180:2);
\draw[red,thick, dashed] (0.83,0) arc (0:94:2.83); 
\draw[red,thick]   (0,2) arc (45:90:2.83) node[left]{$-1+i\sqrt{2}$}   ;  
\node[red] at (-2,2.83){$\bullet$};
\node[red] at (0,2){$\bullet$};
\draw[red,thick, dashed] (-0.83,0) arc (180: 70:2.83); 
\draw[red,thick]   (0,2) arc (135:90:2.83)   ;   

\end{tikzpicture}
\end{center}
\subsection{Twisted embedding}
Put $H'=\left(\begin{smallmatrix}&1\\1&\end{smallmatrix}\right)$ and let
\[
G := g_0\,U(1,1;\ZZ[i])\,g_0^{-1}
\]
be the Cayley transform of $U(1,1;\ZZ[i])$. We define a twisted embedding
\[
\iota\colon G\to Sp(2;\ZZ)
\]
by requiring the following diagram to commute
\[
\xymatrix{U(1,1;\ZZ[i])\ar[r]^{\rho}\ar[d]&Sp(2;\ZZ)\ar[d]\\G\ar[r]^{\iota}&Sp(2;\ZZ)}
\]
where the left vertical arrow is the Cayley transform and the right one is (left) conjugation by the symplectic matrix
\[
t=\left(\begin{smallmatrix}H'&\\ \mathbbm{1}&H'\end{smallmatrix}\right)\in Sp(2;\ZZ);
\]
in terms of generators, we have
\[
\iota\left(\begin{smallmatrix}1&2\\&1\end{smallmatrix}\right)=\left(\begin{smallmatrix}1&&1&\\&1&&1\\&&1&\\&&&1\end{smallmatrix}\right),\ \iota\left(\begin{smallmatrix}&1\\-1&\end{smallmatrix}\right)=\left(\begin{smallmatrix}&-1&1&\\-1&&&1\\-2&&&1\\&-2&1&\end{smallmatrix}\right),\ \iota\left(\tfrac{1+i}{2}\left(\begin{smallmatrix}1&1\\-1&1\end{smallmatrix}\right)\right)=\left(\begin{smallmatrix}1&&&\\-1&&&1\\-1&-1&1&1\\1&-1&&\end{smallmatrix}\right).
\]
It is readily verfied that the induced embedding of symmetric spaces takes the form
\begin{equation}\label{twisted induced embedding}
\iota_{*}\colon\mathbb{H}_1\to\mathbb{S}_{2},\ \tau\mapsto\Omega_{\tau}=\begin{pmatrix}\frac{\tau}{2}&\frac{1}{2}\\\frac{1}{2}&\frac{\tau}{2}\end{pmatrix},
\end{equation}
and the automorphism $J$ is transformed into
\begin{equation}\label{Bolza's linear transformation}
\sigma(t)J\sigma(t)^{-1}=\left(\begin{smallmatrix}0&1&1&0\\-1&0&0&-1\\&&0&1\\&&-1&0\end{smallmatrix}\right).
\end{equation}
Observe that the analog of \eqref{lifted automorphism} reads
\[
( \mathbbm{1},\Omega_{\tau})\cdot \sigma(t)J\sigma(t)^{-1}=\left(\begin{smallmatrix}&1\\-1&\end{smallmatrix}\right)\cdot( \mathbbm{1},\Omega_{\tau}),
\]
and that $\left(\begin{smallmatrix}&1\\-1&\end{smallmatrix}\right)\left(\begin{smallmatrix}i\\1\end{smallmatrix}\right)=-i\left(\begin{smallmatrix}i\\1\end{smallmatrix}\right)$. 
Remarkably, the line spanned by this eigenvector is $G$--invariant; more precisely, noting that for $n=1$ the action \eqref{almost a bundle of Lie algebras} also makes sense for the group $G$, we have the following:
\begin{prop} The map
\[
\mathbb{H}_{1}\times\CC\to\mathbb{S}_{2}\times\CC^{2},\ (\tau,w)\mapsto\left(\Omega_{\tau},\left(\begin{smallmatrix}iw\\w\end{smallmatrix}\right)\right),
\]
is $G$--equivariant. 
\end{prop}
\begin{proof}
This can be shown by a direct computation using the generators.
\end{proof}
\subsection{Automorphic forms}
Recall that the theta group $\Gamma_{\theta}$ is conjugate to the congruence subgroup $\Gamma_{1}(2)\subset SL(2;\ZZ)$, hence its ring of modular forms is also polynomial on two generators; more precisely, we have
\[
M^{\mathbb{C}}_{*}(\Gamma_{\theta})\cong\CC[\delta',\epsilon'],
\]
where $\delta'=\frac{1}{8}(\theta[\begin{smallmatrix}1\\0\end{smallmatrix}](\tau)^{4}-\theta[\begin{smallmatrix}0\\1\end{smallmatrix}](\tau)^{4})$ and $\epsilon'=-\frac{1}{16}\theta[\begin{smallmatrix}1\\0\end{smallmatrix}](\tau)^{4}\theta[\begin{smallmatrix}0\\1\end{smallmatrix}](\tau)^{4}$ are given in terms of Jacobi theta functions. Note that the subscript $*$ refers to a topological grading by \textit{degree}, where we define the degree of a form to be {\em twice} its weight.
\begin{prop} As before, let $G$ be the Cayley transform of $U(1,1;\ZZ[i])$. Then we have
\[
M^{\CC}_{*}(G)\cong\CC[\alpha,\beta],
\]
where $\alpha=-2^{7}(\epsilon'-\frac{1}{2}\delta'^{2})$, $\beta=2^{12}\delta'^{4}$, are  automorphic forms of weight four and eight, respectively.
\end{prop}
\begin{proof}
Clearly, $M^{\CC}_{*}(G)$ is a subring of $M^{\CC}_{*}(\Gamma_{\theta})$. \mbox{In fact, since $\diag(i,i)\in G$,} it follows that $M^{\CC}_{*}(G)$ is concentrated in degrees divisible by eight (i.e.\ in weights divisible by four), hence $M^{\CC}_{*}(G)\subset \CC[\delta'^{2},\epsilon']$. Recall that the element $\frac{1+i}{2}\left(\begin{smallmatrix}1&1\\-1&1\end{smallmatrix}\right)\in G$ induces the map  \eqref{additional transformation}. Making use of the transformation properties (and the duplication formul\ae) for $\theta$--functions, it is straightforward to check that
\[
\delta'({\textstyle{-\frac{2}{\tau-1}-1}})^{2}=+{\textstyle{\frac{(\tau-1)^{4}}{4}}}\delta'(\tau)^{2},\quad \epsilon'({\textstyle{-\frac{2}{\tau-1}-1}})=-{\textstyle{\frac{(\tau-1)^{4}}{4}}}(\epsilon'(\tau)-\delta'(\tau)^{2});
\]
thus, the involution
\[
\star\colon M^{\CC}_{8}(\Gamma_{\theta})\to M^{\CC}_{8}(\Gamma_{\theta}), \quad f(\tau)\mapsto-\tfrac{4}{(\tau-1)^{4}}f({\textstyle{-\frac{2}{\tau-1}-1}})
\]
splits $M^{\CC}_{8}(\Gamma_{\theta})$ into $\pm1$--eigenspaces spanned by $\epsilon'-\frac{1}{2}\delta'^{2}$ and $\delta'^{2}$, respectively. Since $\CC[\delta'^{2},\epsilon']\cong\CC[\delta'^{2},\alpha]$, we conclude that indeed $M^{\CC}_{*}(G)\cong\CC[\delta'^{4},\alpha]\cong\CC[\alpha,\beta]$.
\end{proof}

Each automorphic form admits a Fourier expansion at the (unique) cusp $i\infty$; for all rings $\mathbb{Z} \subseteq R \subseteq \mathbb{C}$ we denote by $M^R_*(G)$ the ring of automorphic forms for $G$ with expansions in $R[\![q^{1/2}]\!]$, where $q=e^{2\pi i\tau}$.

\begin{prop}\label{integral structure}
The subring of automorphic forms expanding integrally at the cusp $i\infty$ is given by
\[
M^{\ZZ}_{*}(G)\cong\ZZ[\alpha,\Delta_{G}],
\]
where $\Delta_{G}=2^{6}\epsilon'(\epsilon'-\delta'^{2})$ is the normalized cusp form.
\end{prop}
\begin{proof}
We have $\alpha^{2}=2^{8}\Delta_{G}+\beta$, hence $\CC[\alpha,\beta]\cong\CC[\alpha,\Delta_{G}]$. Furthermore, from
\[
\delta'=-\textstyle{\frac{1}{8}}+O(q^{1/2}),\quad  \epsilon'=-q^{1/2}+O(q),
\]
and their definitions via Jacobi theta functions it is clear that $8\delta'$ and $\epsilon'$ have integral Fourier expansions at the cusp. Therefore, the same is true for 
\[
\alpha=-2^{7}(\epsilon'-\textstyle{\frac{1}{2}}\delta'^{2})=1+O(q^{1/2}),\quad \Delta_{G}=q^{1/2}+O(q),
\]
and now the claim follows via an easy induction argument.
\end{proof}

\begin{prop}\label{zeroes}
We have
$\alpha(1+i\sqrt{2})=0$ and $\beta(i)=0$. 
\end{prop}
\begin{proof} The element $1+i\sqrt{2}\in\mathbb{H}_{1}$ is fixed under the fractional linear action of $\frac{1+i}{2}\left(\begin{smallmatrix}-1&3\\-1&1\end{smallmatrix}\right)\in G$. Since
\[
\left(\tfrac{1+i}{2}(-(1+i\sqrt{2})+1)\right)^{4}=\left(\tfrac{1-i}{\sqrt{2}}\right)^{4}=-1,
\] 
we conclude that any modular form of weight four (w.r.t.\ $G$) vanishes at $\tau=1+i\sqrt{2}$. For the second statement, note that the transformation properties of the Jacobi $\theta$--functions under $\tau\mapsto-1/\tau$ imply that $\delta'(i)=0$, hence in particular $\beta(i)=0$.
\end{proof}

\begin{prop} 	\label{ka achtel}
Let $f\colon\mathbb{H}_1\to\CC$ be a meromorphic function satisfying  $f(\tfrac{a\tau+b}{c\tau+d})=(c\tau+d)^{k}f(\tau)$ for all $\left(\begin{smallmatrix}a&b\\c&d\end{smallmatrix}\right)\in G$. Then the orders of its poles and zeros satisfy
\[
\textrm{ord}(f,i\infty)+\textstyle{\frac{1}{2}}\textrm{ord}(f,1+i\sqrt{2})+\textstyle{\frac{1}{4}}\textrm{ord}(f,i)+\sum_{p}\textrm{ord}(f,p)=\textstyle{\frac{k}{8}},
\]
where the sum is taken over a set of representatives mod $G$.
\end{prop}
\begin{proof} This can be shown in complete analogy to the case of $SL(2;\ZZ)$ (cf., e.g., \cite{Freitag:2009aa}), i.e., by computing the integral of $\int_{C}f'(\tau)/f(\tau)d\tau$ for a suitable contour $C$. \end{proof}

\begin{prop}\label{compactification}
The function
\[
j_{G}=\frac{\beta}{4\alpha^{2}}\colon(\mathbb{H}_1 \cup \mathbb{Q}P^1)/G\to\CC\cup\{\infty\}\cong\CC P^{1}
\]
is a bijection.
\end{prop}
\begin{proof} The function $j_{G}$ has weight zero and a pole in $1+i\sqrt{2}$. This remains true for any translate $j_{G}-w$, where $w\in\CC$, hence $j_{G}-w$ has a unique zero in $(\mathbb{H}_1 \cup \mathbb{Q}P^1)/G$.
\end{proof}


The following proposition tells us that there is only a single point in the ball quotient where the corresponding abelian two-fold is not a Jacobian. It is the product of two Jacobians, though, as product of two elliptic curves with complex multiplication by the Gaussian integers.

\begin{prop} If a period matrix in the image of $\iota_*$ is diagonalizable by a symplectic transformation, then it is $G$-equivalent to the base point $\iota_*(i)$.
\end{prop}

\begin{proof}
If $\Omega$ is diagonalizable by the fractional linear action of an element in $Sp(2;\ZZ)$, then one of the ten even theta constants vanishes. Now let $\Pi(\Omega)$ be the product of all ten theta constants raised to the eighth power. Then the transformation law of theta functions implies
\[
\Pi\big((A\Omega+B)(C\Omega+D)^{-1}\big)=(\det(C\Omega+D))^{40}\Pi(\Omega)
\]
for all $\left(\begin{smallmatrix}A&B\\C&D\end{smallmatrix}\right)\in Sp(2;\ZZ)$.
Looking at the generators, we deduce that the pullback $\iota^{*}\Pi$ is an automorphic form of weight 80 for $G$. Moreover, the set of the ten even characteristics decomposes into three $G$--orbits, consisting of two, four, and four elements, respectively, hence $\iota^{*}\Pi$ decomposes into a product of automorphic forms of weight 16, 32, and 32, respectively. Computing the first few terms of the Fourier expansions at the cusp, we conclude that $\iota^{*}\Pi$ is proportional to $\beta^{2}\cdot\Delta_{G}^{4}\cdot\Delta_{G}^{4}$; making use of Proposition \ref{ka achtel}, this establishes the claim.
\end{proof}
\section{TAF via hyperelliptic curves}	\label{Section3}
In this section we introduce the family of hyperelliptic curves of genus two we use to get access to the Shimura variety. We use our knowledge about the ring of automorphic forms to give a parametrization of the family in terms of the generators of this ring. A careful analysis of integrality of formal groups then puts us into the position to define a $p$-local genus taking values in this ring. A different choice of local parameter will give a neat description in terms of Legendre polynomials.  Finally, we verify Landweber's criterion.

\subsection{A family of hyperelliptic curves of genus two}
Abelian varieties, even two-folds, are in general difficult to get hands on in terms of defining equations. Thus our approach is to get access to them via Jacobians of curves. Recall that a smooth algebraic curve $C$ of genus $n$ has an associated $n$-dimensional principally polarized abelian variety $J(C)={\rm Pic}^0(C)$, called its Jacobian variety.

Furthermore, recall that any smooth curve of genus two is hyperelliptic, i.e. it is a $2$-fold covering of $\mathbb{P}^1$.
Over $\CC$, any hyperelliptic curve of genus two admits an affine model of the form
\[
 C': y^{2}=p_{5}(x),
\]
where $p_{5}$ is a monic polynomial of degree five with no repeated roots.  While the affine model is an easy way to describe a hyperelliptic curve, its Zariski closure in $\mathbb{P}^2$ has a unique singular point at infinity. It is thus standard practice to model such curves by glueing together two non-singular affine curves $$C= C' \cup C'',$$ one given by the affine model $C'$ above and the other 
$$ C'': v^2 = u^{6}p_5({u}^{-1})$$
by changing coordinates via 
$x=u^{-1}$ and $y=vu^{-3}$. 
In the latter model the points at infinity correspond to points where $u=0$.

\begin{rmk}
An explicit projective embedding of the corresponding Jacobian can be found in \cite{Grant:1990aa}.
\end{rmk}

Each such curve admits the standard hyperelliptic involution $$(x,y)\mapsto(x,-y).$$ The classification of curves of genus two with extra automorphisms dates back to Bolza; as a minor modification of  \cite[\S11, case III]{Bolza:1887aa}, let us introduce the family
\begin{equation}\label{gaussian family}
C'\colon y^{2}=x(x^{4}-2\alpha x^{2}+\beta),
\end{equation}
where we require $\beta(\alpha^{2}-\beta)\neq0$ to ensure smoothness. Each member of this family admits an automorphism of order four, viz.\  $(x,y)\mapsto(-x,iy)$. 

In fact, Bolza exploits this symmetry to deduce a linear transformation on the extended period matrix of the curve, matching the one given by \eqref{Bolza's linear transformation}, and uses this to show that the normalized period matrix $(\tau_{ij})$ satisfies $\tau_{11}=\tau_{22}$, $\tau_{12}=\tfrac{1}{2}$; thus, $(\tau_{ij})$ is in the image of $\iota_{*}$.

\begin{rmk} To exhibit an additional symmetry of these curves, choose $\gamma$ such that $\gamma^{4}=\beta$ and consider the map $(x,y)\mapsto(\gamma^{2}/x,-\gamma^{3}y/x^{3})$, which extends to an involution interchanging the origin and the infinite point.

\end{rmk}

\begin{prop}
Evaluating the automorphic forms $\alpha$ and $\beta$ on the upper half plane, induces a bijection between the points of $\mathbb{H}_{1}/G\setminus\{[i]\}$ and the set of $\CC$--isomorphism classes of curves of genus two whose automorphism group contains the dihedral group of order eight.

\end{prop}
\begin{proof}
It follows from \cite[Proposition 2.1]{Cardona:2007aa} that such a curve is either $\CC$--isomorphic to Bolza's surface, $Y^{2}=X^{5}+X$, or to a unique member of the family $Y^{2}=X^{5}+X^{3}+jX$, where $j\in\CC\setminus\{0,\tfrac{1}{4}\}$.  Clearly, we recover Bolza's surface from \eqref{gaussian family} via $\alpha=0$, i.e.\ $[\tau] = [1+i\sqrt{2}]$, which is the only zero of $\alpha$ due to Proposition \ref{ka achtel}.  If on the other hand $[\tau] \neq [1+i\sqrt{2}]$, we may rescale \eqref{gaussian family} to obtain 
\[
Y^{2}=X^{5}+X^{3}+\frac{\beta(\tau)}{4\alpha(\tau)^{2}}X,
\]
thereby recovering the aforementioned family due to Proposition \ref{compactification}.
\end{proof}

\subsection{The $p$-complete viewpoint}

Note that the curve, its Jacobian, thus all differentials and in particular the splitting of the Lie algebra of the Jacobian are defined over $$R= \mathbb{Z}[i][\tfrac{1}{2}][ \alpha, \beta, (\beta(\alpha^2-\beta))^{-1}].$$ Inverting ${2}$ and  $ \beta(\alpha^2-\beta)$ ensures smoothness of the curve, whereas adjoining the imaginary unit is necessary when considering the $\mathbb{Z}[i]$-actions.

If we complete $R$ at a split prime $p=1\!\!\mod 4$, then $R^\wedge_p$ splits into two copies of $\mathbb{Z}[\alpha, \beta, (\beta(\alpha^2-\beta))^{-1}]^\wedge_p$. Using the fact that the induced action of $\mathbb{Z}[i]$ on the formal group $\mathbb{G}_{J(C)}$ of the Jacobian factors over $R^\wedge_p$, i.e.\ turns ${\rm End}(\mathbb{G}_{J(C)})$ into an $R^\wedge_p$-algebra, splits the formal group into two summands. Each of those $1$-dimensional summands lives over one copy of $\mathbb{Z}[\alpha, \beta,  (\beta(\alpha^2-\beta))^{-1}]^\wedge_p$ and would give rise to a $p$-complete complex oriented cohomology theory via the machinery of \cite{Behrens:2010aa}.

In what follows we use the presentation of the curve to not only refine the above argument to a $p$-local statement, but also give an explicit description of the associated $p$-local genus.
We start analyzing the rational situation.
\subsection{The formal logarithm}
A basis for the holomorphic differentials on a (general) hyperelliptic curve of degree two ($y^{2}=p_{5}(x)$) is given by
\[
\frac{dx}{2y},\ \frac{xdx}{2y}.
\]
Clearly, this basis consists of eigenvectors w.r.t.\ the automorphism of order four on the family \eqref{gaussian family}. If this family  is transformed to the other chart,
\[
C''\colon v^{2}=u(1-2\alpha u^{2}+\beta u^{4}),
\]
we have that $P=(0,0)\in C''$ correponds to the branch point at infinity in the original description $C'$. Via the implicit function theorem, near $P$ the curve can be expressed as the graph over the $v$--line, i.e.\ $v$ locally parametrizes the curve, with
\[		
u(v)=v^{2}+2\alpha v^{6}+(12\alpha^{2}-\beta)v^{10}+O(v^{14}).
\]
Moreover, the differential
\begin{equation}		\label{differential}
\frac{du}{2v}=\frac{dv}{1-6\alpha u^{2}+5\beta u^{4}}
\end{equation}
does not vanish in $P$, while the other differential, $udu/2v$, has a zero in $P$. We define a rational genus
using the parametrization given above, i.e.\ we define
\[
\varphi\colon MU_{*}^{\QQ} \to M_{*}^{\QQ}(G)
\]
to have the formal logarithm
\[
\log_{\varphi}(v)=\int\frac{dv}{1-6\alpha u^{2}+5\beta u^{4}}=v+\frac{6}{5}\alpha v^{5}+O(v^{9}).
\]
We have the following integrality statement:
\begin{klem}
For each prime $p\equiv1\!\mod4$, the image of $BP_*$ under $\varphi$ is contained in the subring of $p$-local automorphic forms:
$$ \varphi(BP_*) \subset M^{\mathbb{Z}_{(p)}}_*(G)  \cong\mathbb{Z}_{(p)}[\alpha, \beta].$$
\end{klem}

\begin{proof}
Having a rational genus and a formal group over $\mathbb{Z}[\alpha, \beta, (\beta(\alpha^2-\beta))^{-1}]^\wedge_p$, we prove the claim making use of an arithmetic fracture square
\begin{center}
\begin{tikzcd}
\mathbb{Z}_{(p)} \ar{d} \ar{r} & \mathbb{Z}_p \ar{d} \\ \mathbb{Q} \ar{r} & \mathbb{Q}_p.
\end{tikzcd}
\end{center}
Choosing a local parameter $v$, as above, at the point at infinity $P$ on the curve, the parametrization $(u(v), v)$ allows us to pull back the $1$-dimensional summand of $\mathbb{G}_{J(C)}$ to a formal neighborhood of the point $P$ in $C''$. In particular, the pullback of the formal differential of the $1$-dimensional summand agrees with the coordinate expression of \eqref{differential} on the nose.
\end{proof}

\subsection{Reparametrization}
We are going to give a more convenient genus taking values in the ring of automorphic forms on $U(1,1;\ZZ[i])$. To this end, recall that the (homogeneous) Legendre polynomials $P_k$ are given via the generating function
\[
(1-2\alpha u+\beta u^{2})^{-1/2}=\sum_{k\geq0}P_{k}(\alpha,\beta)u^{k}.
\]
Now define a rational genus
\[
\varphi^{L}\colon MU_{*}^{\QQ}\to M_{*}^{\QQ}(G)
\]
by means of 
\begin{align*} \log_{\varphi^{L}}'(x)
= \sum_{k=0}^\infty \varphi^{L}[\CC P^{k}]\, x^{k} 
= \sum_{k=0}^\infty P_{k}(\alpha,\beta)\, x^{4k} .
\end{align*} Hence  $\varphi^{L}[\CC P^{l}]=P_{l/4}(\alpha,\beta)$ if $4|l$, and zero otherwise. Then we have the following result:

\begin{thm*} 	\label{Thm}
For each prime $p\equiv1\!\mod4$,  the image of $BP_*$ under $\varphi^L$ is contained in the subring of $p$-local automorphic forms:
$$ \varphi^L(BP_*) \subset M^{\mathbb{Z}_{(p)}}_*(G)  \cong\mathbb{Z}_{(p)}[\alpha, \beta].$$

\end{thm*}
\begin{proof} We choose a different parametrization of $C''$  near $P=(0,0)$, viz.\
\[
u(t)=t^{2},\ v(t)=t(1-2\alpha t^{4}+\beta t^{8})^{1/2}\in\ZZ[\tfrac{1}{2}][\alpha,\beta][\![t]\!];
\]
then we have
\[
\frac{du(t)}{2v(t)}=\frac{dt}{\sqrt{1-2\alpha t^{4}+\beta t^{8}}},
\]
which is the formal differential corresponding to $\varphi^{L}$. As shown before, using the parametrization of the curve $C''$  near $P=(0,0)$ given by $(u(v),v)$, the holomorphic differential $\frac{du}{2v}=\frac{dv}{1-6\alpha u^{2}+5\beta u^{4}}$ gives rise to a one-dimensional formal group law over $M^{\ZZ_{(p)}}_{*}(G)$. But as $v(t)\in\ZZ[\tfrac{1}{2}][\alpha,\beta][\![t]\!]$ is invertible w.r.t.\ composition, $t=t(v)=v+\alpha v^5+O(v^9)$, both descriptions differ only by a reparametrization defined over $\ZZ[\tfrac{1}{2}][\alpha,\beta]\subset M^{\ZZ_{(p)}}_{*}(G)$. Put differently, we have established an isomorphism of commutative 1--dimensional formal group laws over $M^{\ZZ_{(p)}}_{*}(G)$.
\end{proof}

\begin{rmk}
For $\beta=0$, the group law associated to $\varphi^{L}$ reduces to (a variant of) Euler's formal group law $F_{E}$,
\[
F_{E}(x,y)=\frac{x\sqrt{1-2\alpha y^{4}}+y\sqrt{1-2\alpha x^{4}}}{1+2\alpha x^{2}y^{2}}.
\]
This is to be expected, since in this case we are dealing with the group law of a single factor of the product ppav, and this factor is an elliptic curve with multiplication by $i$.

From a different point of view, although for $\beta=0$ the curve defined by the equation \eqref{gaussian family} becomes singular, the projective closure of the (smooth) affine piece $C''$ becomes an elliptic curve with complex multiplication by the Gaussian integers.
\end{rmk}
\subsection{$p$-local $TAF$ over the odd Gaussian lattice}
Having established integrality, we can give explicit descriptions of $p$--local TAF theories. For a fixed prime $p\equiv1\mod4$, let $v_{k}$ denote the image of the $k^{th}$ Hazewinkel generator  \cite[Appendix A2]{Ravenel:2004xh} under $\varphi^{L}$; in particular, we have
\[
v_{1}=P_{\tfrac{p-1}{4}}(\alpha,\beta),\quad v_{2}=\tfrac{1}{p}\big(P_{\tfrac{p^{2}-1}{4}}(\alpha,\beta)-P_{\tfrac{p-1}{4}}(\alpha,\beta)^{p+1}\big).
\]
Recall that $ \alpha^2 - \beta = 2^8 \Delta_G$ is proportional to the normalized cusp form.
\begin{cor}\label{TAF at p=5}
Let $p=5$. Then $\varphi^{L}$ gives $M_{*}^{\ZZ_{(5)}}(G)[\Delta_{G}^{-1}]$ the structure of a Landweber exact $BP_*$-algebra of height two. In particular, the functors
\[
TAF^{U(1,1;\ZZ[i])}_{(5),*}(\cdot)=BP_{*}(\cdot)\otimes_{\varphi^{L}}M^{\ZZ_{(5)}}_{*}(G)[\Delta_{G}^{-1}]
\]
define a homology theory.
\end{cor}
\begin{proof} At the prime 5, the images of the first and second Hazewinkel generator under $\varphi^{L}$ are given by $v_{1}=P_{1}(\alpha,\beta)$ and $v_{2}=(P_{6}(\alpha,\beta)-P_{1}(\alpha,\beta)^{6})/5$. Since $P_{1}(x,1)=x$, $P_{6}(x,1)=\frac{1}{16}(231x^{6}-315x^{4}+105x^{2}-5)$, it follows that $v_{2}\equiv-\beta^{3}\equiv(\alpha^{2}-\beta)^{3}\equiv\Delta_{G}^{3}$ mod $(5,v_{1})$.
\end{proof}

\begin{rmk}
At the prime 5, height two is attained at the ``Bolza point'' $\tau=1+\sqrt{-2}$. Up to $G$--equivalence, this is the only zero of the automorphic form $v_{1}=\alpha$, while $v_{2}$ is non-trivial at this point. If, on the other hand, $\alpha$ and $\beta$ are viewed as indeterminates, these cannot vanish simultaneously, as ensured by invertibility of $\Delta_{G}$; thus height two can only be attained for $\alpha=0$.
\end{rmk}

A tedious but straightforward check reveals that Corollary \ref{TAF at p=5} holds verbatim at the prime 13. We close with a slightly weaker statement valid for an infinite number of primes:

\begin{cor}	\label{CorWeak}
Let $p\equiv5\mod8$, then $v_{2}$ is non-zero mod $(p,v_{1})$. In particular, $M_{*}^{\ZZ_{(p)}}(G)[v_{2}^{-1}]$ admits the structure of a Landweber exact algebra of height two.
\end{cor}
\begin{proof} Note that $\tfrac{p-1}{4}$ is odd for $p\equiv5\mod8$, hence $\alpha|v_{1}$. Making use of the binomial series
\[
(1+\beta t^{8})^{-1/2}=\sum_{n\geq0}\left(\begin{smallmatrix}-1/2\\n\end{smallmatrix}\right)\left(\beta t^{8}\right)^{n}=\sum_{n\geq0}\left(\begin{smallmatrix}2n\\n\end{smallmatrix}\right)\left(-\tfrac{\beta}{4}\right)^{n}t^{8n},
\]
we deduce
\[
pv_{2}\equiv\left(\begin{smallmatrix}\tfrac{p^{2}-1}{4}\\\tfrac{p^{2}-1}{8}\end{smallmatrix}\right)\left(-\tfrac{\beta}{4}\right)^{\tfrac{p^{2}-1}{8}}\mod(\alpha).
\]
Thus, computing the $p$--valuation of the binomial coefficient,
\[
\nu_{p}\left(\begin{smallmatrix}\tfrac{p^{2}-1}{4}\\\tfrac{p^{2}-1}{8}\end{smallmatrix}\right)=\tfrac{p-1}{4}-2\left(\tfrac{p-1}{8}-\tfrac{1}{2}\right)=1,
\]
we conclude that $v_{2}$ is non-zero modulo the ideal $(p,\alpha)$, which in turn implies that $v_{2}$ is non-zero mod $(p,v_{1})$.
\end{proof}

\bibliography{refbib_edited}
\end{document}